\documentclass[12pt]{article}

\usepackage{hyperref}
\usepackage{amsfonts}
\newtheorem{theorem}{Theorem}

\newtheorem{proposition}[theorem]{Proposition}
\newtheorem{lemma}[theorem]{Lemma}
\newtheorem{corollary}[theorem]{Corollary}
\newtheorem{remark}[theorem]{Remark}
\newtheorem{remarks}[theorem]{Remarks}
\newtheorem{example}[theorem]{Example}

\newcommand{\R}{\mathbb{R}}
\newcommand{\C}{\mathbb{C}}

\newcommand{\Sf}{\mathbb{S}}
\newcommand{\grad}{\mbox{grad }}
\newcommand{\spa}{\mbox{span}}

\newcommand{\rank}{\mbox{rank }}

\newcommand{\e}{\epsilon}

\newcommand{\po}{{\hspace*{-1ex}}{\bf .  }}
\newcommand{\ii}{isometric immersion }
\newcommand{\iis}{isometric immersions }

\def\<{\langle}

\def\>{\rangle}
\def\a{\alpha}

\def\d{\partial}

\def\g2{\Gamma^2 }
\def\bea{\begin{eqnarray*} }
\def\eea{\end{eqnarray*} }
\def\be{\begin{equation} }
\def\ee{\end{equation} }

\def\nape{\nabla^\prime}

\def\proof{\noindent{\it Proof: }}
\def\qed{\ifhmode\unskip\nobreak\fi\ifmmode\ifinner\else\hskip5 pt \fi\fi
\hbox{\hskip5 pt \vrule width4 pt  height6 pt depth1.5 pt \hskip 1pt }}
\begin{document}

\title{Euclidean hypersurfaces  with genuine  deformations in codimension two}
\author {M. Dajczer,
  L. Florit\,\,
and\,\, R. Tojeiro
}
\date{}
\maketitle
\

\section{Introduction}

Isometrically  deformable Euclidean hypersurfaces were classified at the
beginning of last century by  Sbrana \cite{sb} and Cartan \cite{ca1}.
 Apart
from flat hypersurfaces, they have rank two, that is,  exactly two nonzero principal
curvatures everywhere, and  are divided into four distinct classes. The two less interesting ones 
consist of ruled  and surface-like
 hypersurfaces. Hypersurfaces in the second class are products with
 Euclidean factors of surfaces in $\R^3$ and cones over surfaces in
 the sphere $\Sf^3$.
The main part of Sbrana's classification is a description of the remaining  classes
in terms  of what is now called  the Gauss parametrization (see \cite{dg1}).
The latter allows to recover a hypersurface of constant rank by means of its
Gauss image and its support function.

The Sbrana-Cartan theory was extended to hypersurfaces of the
sphere and hyperbolic space in \cite{dft}. Moreover, 
 some questions left over
in the works of Sbrana and Cartan, as the possibility of smoothly attaching
different types of hypersurfaces in the Sbrana-Cartan classification
and the existence of hypersurfaces that admit a unique deformation,
were also settled in \cite{dft}.

  The proper setting for attempting to extend the Sbrana-Cartan
theory to submanifolds of higher codimension was developed in \cite{df1}
by means of the concept of {\em genuine\/} isometric  deformations of a submanifold.
An \ii $g\colon\,M^n\to\R^{n+q}$  of an
$n$-dimensional Riemannian manifold $M^n$ into Euclidean space with codimension $q$ is
said to be a {\it genuine isometric deformation\/} of a given \ii
$f\colon\, M^n\to\R^{n+p}$ if $f$ and $g$ are nowhere (i.e.,
on no open subset of $M^n$) compositions, $f=F\circ j$ and $g=G\circ j$, of an isometric
embedding \mbox{$j\colon\,
M^n\hookrightarrow N^{n+r}$}  into a Riemannian manifold $N^{n+r}$
with $r>0$ and \iis $F\colon\,N^{n+r}\to\R^{n+p}$ and $G\colon\,N^{n+r}\to\R^{n+q}$:

\bigskip

\begin{picture}(110,84)\nonumber
\put(115,31){$M^n$}
\put(169,31){$N^{n+r}$}
\put(202,62){$\R^{n+p}$}
\put(202,0){$\R^{n+q}$}
\put(155,59){$f$}
\put(155,4){${g}$}
\put(200,46){${}_F$}
\put(199,26){${}_{G}$}
\put(153,40.5){${}_j$}
\put(188,42){\vector(1,1){16}}
\put(188,28){\vector(1,-1){16}}
\put(137,44){\vector(3,1){60}}
\put(137,26){\vector(3,-1){60}}
\put(145,34){\vector(1,0){21}}
\put(145,36){\oval(7,4)[l]}
\end{picture}
\bigskip\medskip

\noindent An \ii $f\colon\,M^n\to\R^{n+p}$ with the property that  no open subset
$U\subset M^n$
admits a genuine isometric deformation in $\R^{n+q}$ for a fixed  $q> 0$ is said to be
{\it genuinely rigid} in $\R^{n+q}$.

The isometric deformation problem for Euclidean submanifolds is then to classify all 
isometric immersions $f\colon\,M^n\to\R^{n+p}$ that admit genuine isometric deformations 
in $\R^{n+q}$, for each fixed $q>0$. This is a rather difficult problem that has only been solved
in  particular cases. It was shown in \cite{df1} that  if an
\ii $f\colon\, M^n\to\R^{n+p}$ and a genuine isometric deformation
$\hat f\colon\,M^n\to\R^{n+q}$ of it have sufficiently low codimensions
then they are {\it mutually (isometrically) ruled}, that is, $M^n$
carries an integrable \mbox{$d$-dimensional} distribution $D^d\subset TM$
whose leaves are mapped diffeomorphically by $f$ and $\hat f$ onto open
subsets of affine subspaces of $\R^{n+p}$ and $\R^{n+q}$,
respectively. A sharp estimate on the dimension
$d$ of the rulings was also given, and it was proved that the normal connections and second
fundamental forms of $f$ and $\hat f$ satisfy strong additional
relations. However, these necessary conditions for the existence of genuine isometric deformations are 
far from being sufficient, as already shown by the classical Sbrana-Cartan theory for the special case $p=1=q$.

The case $p=2=q$ is the only other one in which some partial results on this problem have been obtained (see \cite{df2} and \cite{dm}), and 
just for the simplest class of submanifolds that satisfy the necessary conditions in \cite{df1} referred to above, namely,  submanifolds of rank two. These can be divided
into three distinct subclasses, called {\em parabolic\/}, {\em elliptic\/}  and {\em
hyperbolic\/}. 

Parabolic submanifolds $g\colon\,M^n\to\R^{n+2}$  that are ruled 
admit genuine isometric deformations of a special type (see \cite{df2}), whereas those that are neither ruled nor surface-like were shown in \cite{dm} to be not only genuinely rigid, but even isometrically rigid in the standard sense.  On the other hand, it was shown in \cite{df2} that any  elliptic submanifold $g\colon\,M^n\to\R^{n+2}$  is  genuinely rigid in $\R^{n+2}$,
unless $M^n$ admits an isometric immersion into $\R^{n+1}$. However,
it was not clear whether there exists any elliptic submanifold  in $\R^{n+2}$ that 
 can be isometrically immersed into Euclidean space as a hypersurface. This was one of the motivations of the present work, and naturally leads us to consider the problem in the case $p=1$ and $q=2$. 


  Notice that an  \ii $g\colon\,M^n\to\R^{n+q}$  is
a genuine isometric deformation of a hypersurface
$f\colon\, M^n\to\R^{n+1}$ if $g$ is nowhere a composition $g=H\circ f$, where
$H\colon\,V\subset \R^{n+1}\to\R^{n+q}$ is an \ii
of an open subset $V$ containing $f(M)$.
In terms of this concept, the main result of \cite{dt1} states that a necessary condition
for a hypersurface  $f$ in $\R^{n+1}$ to admit a genuine isometric deformation
in $\R^{n+q}$, $2\leq q\leq n-2$, is that $\rank f\leq q+1$.
It follows that  $f$ must have rank at most $3$ in order to admit a genuine isometric
deformation in  $\R^{n+2}$.

In this paper, we give a complete local  description, in terms of the Gauss
parametrization,
of rank two Euclidean hypersurfaces of dimension $n\geq 3$
that admit a genuine isometric deformation in $\R^{n+2}$. They are divided into
three distinct classes. Again, ruled and surface-like hypersurfaces form the less
interesting ones. On the other hand, the characterization of the Gauss images
of hypersurfaces in the  remaining class turns out to be significantly more
involved than that in the Sbrana-Cartan theory. 
Even though our result shows that such hypersurfaces 
are rather special, it  allows us to easily construct a large family of explicit examples. 

\section[The  result]{The result}

In order to state our  result, we first recall the Gauss parametrization of
Euclidean hypersurfaces with constant rank.
\medskip

Let $f\colon\,M^n\to\R^{n+1}$ be an oriented hypersurface  such that its Gauss map
$N\colon\, M^n\to\Sf^n$
has constant rank $k$.  Let  $\pi\colon\,M^n\to L^{k}$
denote the projection onto the   quotient space $L^k:=M^n/\Delta$ of totally geodesic
leaves of the
relative nullity distribution of $f$ defined by
$$
\Delta(x)=\ker A(x).
$$
Here $A=A_N$ stands for the shape operator of $f$, given by $AX=- \nabla_XN$
for any $X\in TM$. The
quotient $L^k$ is a differentiable manifold, except possibly for the
Hausdorff condition, but this is always satisfied if we work  locally.

Since the Gauss map is constant along the leaves of $\Delta$,
it induces an immersion
$h\colon\,L^{k}\to\Sf^n$ given by
$
h\circ\pi = N,
$
which we will also call the Gauss map of $f$.
Accordingly, the  support function $\tilde\gamma=\<f, N\>$ of $f$
gives rise to $\gamma\in C^\infty(L)$ defined by $\gamma\circ\pi=\tilde
\gamma$.

The  \emph{Gauss parametrization} given in \cite{dg1} allows to locally recover $f$  in
terms
of the pair $(h,\gamma)$. Namely, there exists locally a  diffeomorphism
$\Phi\colon\,V\subset T_h^\perp L\to M^n$  of an open neighborhood of the zero
section of  the normal bundle $T_h^\perp L$ of $h$ such that $\pi\circ\Phi=\hat \pi$, where $\hat \pi\colon\,T_h^\perp L\to L$ is the
canonical projection, and 
$$
f\circ \Phi(y,w)= \gamma(y) h(y) + h_*\,\nabla\gamma(y) + w.
$$

Next, we introduce the class of spherical surfaces that appear as Gauss images
of hypersurfaces of rank two admitting genuine isometric deformations in codimension two.
\vspace{1ex}

  Let $h\colon\,L^2\to\Sf^n$ be a surface and let   $ \alpha_h$ denote  its second
fundamental form
with values in its normal bundle. Assume $n\geq 4$ and that the first normal spaces
$N^1_h$
of $h$, that is, the subspaces of the normal spaces spanned by the image of $\alpha_h$, 
have dimension 
two everywhere.
Then,  given a frame $\{X, Y\}$ of $TL$,
there exist  $a, b, c \in \R$ such that
$$
a\alpha_h (X,X) + 2c \alpha_h(X, Y ) + b \alpha_h(Y, Y ) = 0.
$$
We say that $h$ is \emph{elliptic} (resp.,  \emph{hyperbolic} or \emph{parabolic})  if
$ab- c^2 > 0$
(resp., $< 0$ or $= 0$) everywhere, a condition that is independent of the given frame.
This is easily seen to be equivalent
to the existence of  a tensor $J$ on $L^2$ satisfying $J^2=\e I$, with $\e=-1$, $1$ or $0$
respectively, and
\be\label{eq:aJ}
\alpha_h(JX, Y)=\alpha_h(X, JY).
\ee
 Moreover, in the first two cases the tensor $J$ is unique up to sign. Note that
(\ref{eq:aJ})
is equivalent to requiring  any height function $h^v=\<h,v\>$
of $h$ to satisfy
$$
({\rm Hess}_{\,h^v}+h^v I)J=J^t({\rm Hess}_{\,h^v}+h^v I).
$$
Here ${\rm Hess}$ denotes the Hessian with respect to metric induced by $h$,  as an
endomorphism of
$TL$. We also use the same notation  for the corresponding symmetric bilinear form.

When the first normal space of $h$ has dimension less than two, we still
call $h$ elliptic, hyperbolic or parabolic with respect to such a
tensor $J$ if (\ref{eq:aJ}) is satisfied.


  Now assume that   $h\colon\,L^2\to\Sf^n$ is hyperbolic. Let $(u,v)$
be coordinates around $(0,0)$ whose coordinate vector fields $\{\d_u, \d_v\}$ are
eigenvectors of $J$.
Then, condition  (\ref{eq:aJ}) says that $(u,v)$ are  conjugate
coordinates for $h$, that is,
 $$
\alpha_h(\d_u,\d_v)=0.
$$
 Write
$$
\nabla_{\partial_u}\partial_v=\Gamma^u\partial_u+\Gamma^v\partial_v,
$$
and set $F=\langle\partial_u,\partial_v\rangle$. Then, any height function $h_v$
of $h$ belongs to the kernel of the differential operator
$$
Q(\theta):= {\rm Hess}_{\,\theta}(\partial_u,\partial_v)+F\theta
=\theta_{uv}-\Gamma^u \theta_u - \Gamma^v \theta_v + F\theta.
$$

For each pair of smooth functions $U=U(u)$ and $V=V(v)$   define
$$
\varphi^U=U(u)e^{2\int_0^v\Gamma^u(u,s)ds}\ \  \mbox{and} \ \
\psi^V=V(v)e^{2\int_0^u\Gamma^v(s,v)ds}.\ \ \ \
$$
In other words, the functions $\varphi^U$ and $\psi^V$ satisfy
\be\label{eq:vp}
\varphi^U_v=2\Gamma^u\varphi^U\ \  \mbox{and} \ \ \psi^V_u=2\Gamma^v\psi^V
\ee
as well as the initial conditions
$\varphi^U(u,0)=U(u)$ and $\psi^V(0,v)=V(v)$.

Assume, in addition, that one of the following conditions hold:
\be\label{eq:cond1}
  U,V>0 \,\,\,\,\mbox{or}\,\,\,\, 0<2\varphi^U<-(2\psi^V+1)\,\,\,\,
\mbox{or}\,\,\,\, 0<2\psi^V<-(2\varphi^U+1).
\ee
Define
$$
\rho^{UV}=\sqrt{|2(\varphi^U+\psi^V)+1|}
$$
and
$$
{\cal C}_h=\{(U,V): (\ref{eq:cond1}) \mbox{ holds and }
Q(\rho^{UV})=0\}.
$$

Consider now an elliptic surface  $h\colon\,L^2\to\Sf^n$. If $(u,v)$
are coordinates on $L^2$ such that
$$
\d_z=\frac{1}{2}(\d_u-\d_v)\,\,\,\,\mbox{and}\,\,\,\, \d_{\bar z}=\frac{1}{2}(\d_u+\d_v)
$$
are eigenvectors of the complex linear extension of $J$ to $TL\otimes \C$,
where $\{\d_u, \d_v\}$ is the frame of coordinate vector fields, then
$$
\alpha_h(\d_z, \d_{\bar z})=0.
$$
Define a complex--valued Christoffel symbol $\Gamma$ by
$$
\nabla_{\d_ z}\d_{\bar z}=\Gamma \d_z+ \bar \Gamma \d_{\bar z}
$$
and set $F=\<\partial_z,\partial_{\bar z}\>$, where $\<\,\,, \,\,\>$ also stands
for the complex bilinear extension of the metric induced by $h$.
In this case, the height functions
of $h$ belong to the kernel of the differential operator
$$
Q(\theta): = {\rm Hess}_{\,\theta}(\partial_z,\partial_{\bar z})+F\theta
=\theta_{z\bar z}-\Gamma \theta_z - \bar \Gamma \theta_{\bar z} + F\theta.
$$
For each holomorphic function $\zeta$, let $\varphi^{\zeta}(z,\bar z)$ be  the unique
complex valued function such that
$$
\varphi^{\zeta}_{\bar z}=2\Gamma{\varphi}^{\zeta}\ \ \mbox{and}\ \
\varphi^{\zeta}(z, 0)=\zeta(z).
$$
Assume further that
\be\label{eq:cond2}
\varphi^{\zeta}\neq -1/2\,\,\,\,\mbox{and}\,\,\,\,\,4\mbox{Re}(\varphi^{\zeta}) +1<0.
\ee
In this case, define
$$
\rho^{\zeta}=\sqrt{-(4\mbox{Re}(\varphi^{\zeta}) +1)}
$$
and
$$
{\cal C}_h=\{\zeta\ {\rm holomorphic} : (\ref{eq:cond2}) \mbox{ holds and }
Q(\rho^{\zeta})=0\}.
$$

For the statement of  our main result, we need a few more definitions.
That an isometric immersion $g\colon\,M^n\to\R^{n+2}$ is {\em locally substantial\/}
 means that there exists no open subset $U\subset M^n$ such that $g(U)$ is contained
in an affine hyperplane of $\R^{n+2}$.
A hypersurface  $f\colon\,M^n\to\R^{n+1}$ is   {\it ruled\/} if $M^n$
admits a foliation by leaves of codimension one that are mapped by $f$ into
affine subspaces of $\R^{n+1}$. It is {\em surface like\/} if $f(M)\subset
L^2\times\R^{n-2}$
where $L^2$ is a surface in $\R^3$,
or $f(M)\subset CL^2\times\R^{n-3}$ where $CL^2\subset\R^4$ is the cone over a surface
$L^2\subset\Sf^{\,3}$.
\vspace{1ex}

\begin{theorem}\po\label{thm:main}
Let $f\colon\, M^n\to\R^{n+1}$ be a nowhere surface--like or ruled rank
two hypersurface that admits a locally substantial genuine isometric
deformation $g\colon\,M^n\to\R^{n+2}$. Then, on each connected component of an open dense subset, the Gauss map $h$ of $f$ is either an elliptic or hyperbolic
surface, ${\cal C}_h$ is nonempty, and the support function $\gamma$
of $f$ satisfies $Q(\gamma)=0$.

Conversely, any simply connected hypersurface $f$ whose Gauss map $h$
and support function $\gamma$ satisfy the above conditions admits
genuine isometric deformations in $\R^{n+2}$, which are parametrized by
the set ${\cal C}_h$.
\end{theorem}

%

\begin{remarks}\po {\em $(i)$ Excluding surface--like and ruled hypersurfaces 
has the only purpose of emphasizing the really interesting class of hypersurfaces that
admit substantial genuine isometric
deformations in $\R^{n+2}$.  In fact, it was shown in \cite{df2} that ruled
hypersurfaces in $\R^{n+1}$ admit locally as many genuine isometric deformations in $\R^{n+2}$ as
triples of smooth functions
in one variable, all of them being also ruled with the same rulings. On the other hand,
if $f$ is surface--like then  genuine deformations of $f$ in $\R^{n+2}$  are
given by  genuine deformations of $L^2$ in either $\R^4$ or $\Sf^4$ (see Remark
\ref{re:surfacelike}).\vspace{1ex}\\
$(ii)$ For $n\geq 4$ and $h$ as in Theorem \ref{thm:main}, it follows from Theorem $1$ in \cite{dt2} that if  $N_1^h$ is one-dimensional everywhere and $h(L)$ is not contained in a totally geodesic 
 $\Sf^3\subset \Sf^n$, then $h$ carries a relative nullity distribution of rank one. Moreover, it is easily seen that in this case
 $h$ is necessarily hyperbolic. }
\end{remarks}

It follows from Theorem \ref{thm:main} that hypersurfaces 
that are neither  surface--like nor ruled and  admit  locally substantial
genuine isometric deformations in $\R^{n+2}$
are rather special. First, if $n\geq 4$
then their Gauss maps
$h$ must have  first normal spaces with dimension less than or equal to two everywhere, which
is
already a strong restriction for \mbox{$n\geq 5$.} Moreover,  $h$ can
not be parabolic, and even if it is elliptic or hyperbolic the condition that
 the set ${\cal C}_h$  be nonempty occurs only in very special cases. Furthermore, one has the condition  $Q(\gamma)=0$ on the support function. Still, we
show next that our  result can be used to easily construct a
family of hyperbolic hypersurfaces  admitting a large set of
genuine deformations in codimension two. It is very likely that similar examples exist in the elliptic case,
although this remains an open problem.

\begin{example}\po {\em Let us analyze the case where
$h$ has flat normal bundle, that is,
$(u,v)$ are real orthogonal conjugate coordinates. Setting
$E=\<\d_u, \d_u\>$ and $G=\<\d_v, \d_v\>$, it follows from $F=\<\d_u,\d_v\>=0$ that
$$
2\Gamma^u=E_v/E\,\,\,\,\,\mbox{and}\,\,\,\,2\Gamma^v=G_u/G.
$$
Hence $\varphi^U=UE$ and $\varphi^V=VG$, after replacing the smooth functions $U(u)$ and
$V(v)$
if necessary. Moreover,
$$
Q(\theta)=\theta_{uv}-\frac{E_v}{2E}\theta_u-\frac{G_u}{2G}\theta_v.
$$
Assume further that $(u,v)$ are also isothermic coordinates, say,  $E=e^{2\lambda}=G$, so
that $h$ is an
isothermic surface. Then $Q(\theta)=\theta_{uv}-\lambda_v\theta_u-\lambda_u\theta_v$.
Hence, we have ${\cal C}_h\neq \emptyset$ if and only if there exist smooth functions
$U(u)$ and
$V(v)$ such that $\rho:=\sqrt{(U+V)e^{2\lambda}+1}$ satisfies $Q(\rho)=0$, that is,
$\rho_{uv}=\lambda_v\rho_u+\lambda_u\rho_v.$
The latter equation is equivalent to
$$
4(U+V)(1-e^{2\lambda}(U+V))\lambda_{uv}-4e^{2\lambda}(U+V)^2\lambda_u\lambda_v+2U'\lambda_v+2V'\lambda_u
-e^{2\lambda}U'V'=0.
$$
In addition, now suppose that  $\lambda=\lambda(u)$. Then, these isothermic surfaces are
precisely
warped products of curves in the sense of \cite{no}, parametrized in isothermic
coordinates
(see the main theorem in \cite{no}). Then,
the preceding equation reduces to
$$
V'(2\lambda_u-e^{2\lambda}U')=0.
$$
It follows that ${\cal C}_h$ is the set of pairs of
smooth functions $(U,V)$ satisfying
$$
 V=d\in\R\,\,\,\,\,{\mbox or}\,\,\,\, U=c-e^{-2\lambda}
\,\,\,\mbox{for}\,\,c\in \R,
$$
with the restrictions arising from (\ref{eq:cond1}).
By Theorem \ref{thm:main}, any simply connected hypersurface given in terms
of the Gauss parametrization by a pair $(h,\gamma)$, where $h$ is a surface as above and
$Q(\gamma)=\gamma_{uv}-\lambda_u\gamma_v=0$, i.e., $\gamma=\nu e^\lambda$ for some smooth
function $\nu(v)$, admits genuine isometric deformations in $\R^{n+2}$ which are
parametrized by
${\cal C}_h$.
}\end{example}

\section[Outline of the proof]{Outline of the proof}

The starting point of the proof of Theorem \ref{thm:main} is  the following  consequence
of Proposition $7$ in  \cite{df2}.

\begin{proposition}\po\label{le:0} Let
$f\colon\, M^n\to\R^{n+1}$ be a rank two hypersurface. Then,
any genuine isometric deformation of $f$ in $\R^{n+2}$ has the same
relative nullity distribution as $f$.
\end{proposition}

The remaining of the proof is divided into four main steps. In the first
one, we approach the problem in the light of the fundamental theorem of
submanifolds. We show that the existence of a locally substantial
isometric immersion $g\colon\,M^n\to \R^{n+2}$ with the same relative
nullity distribution as $f$ implies (and is implied by, if $M^n$ is
simply connected) the existence of a pair of tensors and a one-form
satisfying several conditions that arise from the Gauss, Codazzi and Ricci
equations for an isometric immersion into $\R^{n+2}$.

A key role in the proof is played by the splitting tensor of the relative
nullity distribution $\Delta$ of $f$. Recall that the {\em splitting
tensor\/} $C$ of a totally geodesic distribution $\Delta$ on a Riemannian
manifold $M^n$ assigns to each $T\in\Delta$ the endomorphism $C_T$ of
$\Delta^{\perp}$ given by
$$
C_TX = -(\nabla_{X}T)_{\Delta^{\perp}}.
$$
   When $\Delta^\perp$ has rank two, we say that $\Delta$ is  \emph{elliptic},
\emph{hyperbolic} or \emph{parabolic}  if there exists a tensor
$J\colon\,\Delta^\perp\to \Delta^\perp$
satisfying $J^2=\epsilon I$, with $\epsilon=-1$, $1$ or $0$, respectively, such that the
image
$C(\Delta)$ of $C$ at each point of $M^n$ is not spanned by the identity tensor $I$ but is
contained in $\spa\{I,J\}$. Accordingly, we call an Euclidean submanifold of rank two
\emph{elliptic}, \emph{hyperbolic} or \emph{parabolic} if its relative nullity
distribution
is elliptic, hyperbolic or parabolic, respectively.

 \begin{remark}\po{ \em   Let $g\colon\, M^n\to \R^{n+p}$  be an isometric immersion of
rank two of a Riemannian manifold without
flat points. Assume that the first normal spaces of $g$ have dimension two everywhere.
As in the case of surfaces with  first normal bundle of rank two, there exists a tensor
$J\colon\,\Delta^\perp\to \Delta^\perp$, where $\Delta$ is the relative nullity
distribution,
satisfying $J^2=\epsilon I$, with $\epsilon=-1$, $1$ or $0$, and $\alpha(JX, Y)=\alpha(X,
JY)$.
In \cite{df1} the immersion $g$ was called accordingly \emph{elliptic}, \emph{hyperbolic}
or
\emph{parabolic}, respectively, and it was shown that this is equivalent to the preceding
definition.}
\end{remark}

\begin{proposition}\po\label{le:1}
Let $f\colon\, M^n\to\R^{n+1}$ be a rank two hypersurface that is
nowhere surface-like or ruled.
Assume that there exists a locally substantial isometric immersion
$g\colon\,M^n\to\R^{n+2}$
that also has $\Delta$ as its relative nullity distribution. Then, on
an open dense subset of $M^n$, the hypersurface $f$ is either elliptic or hyperbolic and
there exist a unique (up to signs and permutation) pair $(D_1,D_2)$ of tensors in
$\Delta^\perp$
contained in $\spa\,\{I,J\}$
and a unique  one-form $\phi$ on $M^n$ satisfying  the following conditions:
\begin{itemize}
\item[(i)] $\Delta\subset \ker\phi$,
\item[(ii)] $AD_i=D_i^tA$,
\item[(iii)]  $\det D_i=1/2$,
\item[(iv)] $\nabla_T D_i=0=[D_i, C_T]$ for all $T\in \Delta$,
\item[(v)]
$(\nabla_{X}AD_i)Y-(\nabla_{Y}AD_i)X=(-1)^{j}A(\phi(X)D_{j}Y-\phi(Y)D_{j}X), \,\, i\neq
j$,
\item[(vi)]  $d\phi(Z,T)=0$ for all $Z\in TM$ and $T\in \Delta$,
\item[(vii)] $d\phi(X,Y)=\<[AD_1, AD_2]X,Y\>$,
\item[(viii)]  $D_2^2\neq \pm D_1^2$.
\end{itemize}

   Conversely, if $M^n$ is simply connected and  $f$ is an  elliptic or hyperbolic
hypersurface
that carries such a triple $(D_1, D_2, \phi)$, then there exists  a  locally substantial
isometric immersion $g\colon\,M^n\to\R^{n+2}$   with
$\Delta$ as its relative nullity distribution.
Moreover, distinct triples (up to signs and permutation) yield
noncongruent isometric immersions, and conversely.
\end{proposition}

The second step is to show that the problem of finding a triple $(D_1,
D_2, \phi)$ satisfying all conditions in Proposition \ref{le:1} can be
reduced to a similar but easier problem for the Gauss map $h$ of $f$.

\begin{proposition}\po\label{le:2}
Let $f\colon\, M^n\to\R^{n+1}$ be a hypersurface of rank two given in terms of
the Gauss parametrization by a pair  $(h,\gamma)$. Let
$\pi\colon\, M^n\to L^2$ be the projection and $\nape$ the Levi-Civita
connection on $L^2$ for the metric $\<\,\,,\,\,\>'$ induced by $h$.
If $f$ is elliptic (resp., hyperbolic) and  $(D_1, D_2, \phi)$ is a triple on $M^n$
satisfying all conditions in Proposition~\ref{le:1}, then  there exist a unique tensor
$\bar J$ on $L^2$ with $\bar J^2=-I$ (resp., $\bar J^2=I$), a unique  pair of  tensors
$(\bar D_1, \bar D_2)$ in $\spa\,\{I,\bar J\}$, and a
unique  one-form $\bar \phi$ on $L^2$ such that $h$ is elliptic (resp., hyperbolic)
with respect to $\bar J$ and
$$
\bar J\circ \pi_*=\pi_*\circ J,\,\,\,\,\,\,
\bar D_i\circ \pi_*=\pi_*\circ D_i\ \  \mbox{and}\ \
\bar\phi\circ\pi_*=\phi.
$$
Moreover,
\be\label{eq:hessgamma}
({\rm Hess}_{\,\gamma}+ \gamma I)\bar J
=\bar J^t\,({\rm Hess}_{\,\gamma}+ \gamma I)
\ee
and the triple $(\bar D_1, \bar D_2, \bar \phi)$  satisfies:
\begin{itemize}
\item[(a)] $\det \bar D_i=1/2$,
\item[(b)] $(\nabla'_{X}\bar D_i)Y-(\nabla'_{Y}\bar
D_i)X=(-1)^{j}(\bar\phi(X)\bar D_{j}Y-\bar\phi(Y)\bar D_{j}X),\ i\neq j$,
\item[(c)] $d\bar\phi(X,Y)=\<[\bar D_1,\bar D_2]X,Y\>'$,
\item[(d)]  $\bar D_2^2\neq \pm \bar D_1^2$.
\end{itemize}

Conversely,
if $h$ is elliptic (resp., hyperbolic) then $f$ is elliptic (resp., hyperbolic) and
any such triple $(\bar D_1, \bar D_2, \bar \phi)$ on $L^2$ gives
rise to a unique triple $(D_1, D_2, \phi)$ on $M^n$ satisfying all conditions in
Proposition \ref{le:1}.
\end{proposition}

In the third step, we determine under which additional conditions the
isometric deformation $g$ of $f$ in Proposition \ref{le:1} is genuine.

\begin{proposition}\label{le:3}\po
Under the assumptions of Proposition \ref{le:1}, we have that $g$ is a genuine
isometric deformation of $f$ if and only if $\rank ( D_1^2+ D_2^2-I)=2$. This is
always the case if $f$ is elliptic.
\end{proposition}

   The last and crucial step in the proof of Theorem \ref{thm:main} is to characterize the
pairs $(h,\gamma)$,  where $h\colon\,L^2\to \Sf^n$ is an elliptic or hyperbolic surface
and $\gamma\in C^\infty(L)$,
such that $L^2$ carries a triple $(\bar D_1, \bar D_2, \bar \phi)$ satisfying
all conditions in Proposition \ref{le:2} and $\rank (\bar D_1^2+ \bar D_2^2-I)=2$.

\begin{proposition}\po\label{le:4}
Let $h\colon\,L^2\to \Sf^n$ be an elliptic or hyperbolic  surface and let  $\gamma\in
C^\infty(L)$.
Then there exists a triple  $(\bar D_1, \bar D_2, \bar \phi)$
satisfying  all conditions
in Proposition~\ref{le:2} and $\rank (\bar D_1^2+ \bar D_2^2-I)=2$
if and only if  ${\cal C}_h$ is nonempty
and $\gamma$ satisfies   $Q(\gamma)=0$. Distinct triples (up to signs and permutation)
give rise
to distinct elements of  ${\cal C}_h$, and conversely.
\end{proposition}

The proofs of Propositions \ref{le:1} to \ref{le:4} will be provided in the following
sections.
Theorem \ref{thm:main} follows easily by putting them together with
Proposition~\ref{le:0}.
\vspace{2ex}

\noindent {\em Proof of Theorem \ref{thm:main}:} Under the assumptions of Theorem
\ref{thm:main},
it follows from Propositions \ref{le:0},  \ref{le:1} and  \ref{le:3}  that, on an open
dense subset of $M^n$,
the hypersurface $f$ is elliptic or hyperbolic and there exists a unique (up to signs and
permutation)
triple $(D_1, D_2, \phi)$ satisfying all conditions in   Proposition \ref{le:1} as well
as that
in Proposition \ref{le:3}.
Let $f$ be locally given in terms of  the Gauss parametrization by a pair $(h,\gamma)$.
By  Proposition \ref{le:2},
the surface $h$ is elliptic or hyperbolic, respectively, and the triple $(D_1, D_2, \phi)$
projects to  a (unique)  triple  $(\bar D_1, \bar D_2, \bar \phi)$ on $L^2$  satisfying
all conditions in   Proposition \ref{le:2}.
Moreover, $\rank (\bar D_1^2+ \bar D_2^2-I)=\rank (D_1^2+  D_2^2-I)=2$.
We conclude from Proposition~\ref{le:4} that $(\bar D_1, \bar D_2, \bar \phi)$ gives rise
to a unique element of ${\cal C}_h$, and that
$Q(\gamma)=0$.

  Conversely, assume that $f\colon\,M^n\to \R^{n+1}$ is a simply connected hypersurface
given in terms of the Gauss parametrization by a pair $(h,\gamma)$, where
$h$ is an elliptic or hyperbolic surface such that ${\cal C}_h$ is nonempty
and $Q(\gamma)=0$. By  Proposition \ref{le:4}, each  element of ${\cal C}_h$ gives rise to
a unique triple $(\bar D_1, \bar D_2, \bar \phi)$ on $L^2$  satisfying  all conditions in
Proposition \ref{le:2} and $\rank (\bar D_1^2+ \bar D_2^2-I)=2$.  Then, it follows
from Proposition \ref{le:2} that $f$ is elliptic or hyperbolic, respectively,  and that
$(\bar D_1, \bar D_2, \bar \phi)$ can be lifted  to a unique triple $(D_1, D_2, \phi)$ on
$M^n$
satisfying all conditions in   Proposition \ref{le:1}. Moreover, $\rank (D_1^2+ 
D_2^2-I)=2$.
By Proposition \ref{le:1}, such triple
yields a unique (up to rigid motions of $\R^{n+2}$)  isometric immersion $g\colon\,M^n
\to \R^{n+2}$
sharing with $f$ the same relative nullity distribution.  Proposition \ref{le:3} implies
that $g$
is a genuine isometric deformation of $f$ in $\R^{n+2}$.

Finally,   we also have from Propositions \ref{le:0} to \ref{le:4}  that (congruence
classes of) genuine
isometric deformations of $f$ are in one-to-one correspondence with triples $(\bar D_1,
\bar D_2, \bar \phi)$
on $L^2$  satisfying  all conditions in   Proposition~\ref{le:2} and
$\rank (\bar D_1^2+ \bar D_2^2-I)=2$ (up to signs and permutation),
and these are in turn in  one-to-one
correspondence with elements of ${\cal C}_h$. \qed

\section{Projectable tensors and one-forms}

In this section, we  establish some  facts that will be needed in the proof of
Proposition \ref{le:2}
and also have interest on their own.\vspace{1ex}

Given a submersion $\pi\colon\,M\to L$, a vector field $X$ on $M$ is said  to be {\em
projectable\/}
if it is $\pi$-related to a vector field on $L$, that is, if  there exists a vector
field $\bar X$
on $L$ such that $\pi_*X=\bar X\circ \pi$.

\begin{proposition}\po\label{prop:proj}
Let $\Delta$ be an integrable distribution on a differentiable
manifold $M$,  let $L=M/\Delta$ be the (local) quotient space of leaves of $\Delta$ and
let
$\pi\colon\,M\to L$ be the projection. Then, a vector field $X$ on $M$ is projectable if
and only if
$[X,T]\in \Delta$ for any $T\in \Delta$.
\end{proposition}
\proof If $\pi_*X=Z\circ \pi$ and $T\in \Delta$, then
$\pi_*[X,T]=[\pi_*X,\pi_*T]=[Z,0]\circ \pi=0$,
hence $[X,T]\in \Delta$. For the converse, in order to prove that $X$ is projectable we
must show that,
for each leaf $F$ of $\Delta$, the map $\psi\colon\, F\to T_qL$, $q=\pi(F)$, given by
$\psi(p)=\pi_*(p)X_p$,
is constant. Given $p\in F$ and $v\in T_pF$, choose $T\in \Delta$ with $T(p)=v$ and let
$g_t$ be
the flow of $T$.  By the assumption and since $\pi \circ g_t = \pi$,  we have
$$
\begin{array}{l}{\displaystyle 0 = \pi_*[X,T](p) = \lim_{t\mapsto 0} \frac{1}{t}(\pi_*X(
g_t(p))
- \pi_* {g_t}_* X(p))}\vspace{1ex}\\
\hspace*{2ex} {\displaystyle = \lim_{t\mapsto 0} \frac{1}{t}(\pi_*X( g_t(p)) - \pi_* X(p))
=\psi_*(p)v. \qed} \end{array}
$$

If we apply Proposition \ref{prop:proj} to  a totally geodesic distribution $\Delta$ on a
Riemannian
manifold $M$, the conclusion can be expressed in terms of its  splitting tensor $C$.

\begin{corollary}\po\label{cor:proj}
Let $\Delta$ be a totally geodesic distribution on a Riemannian
manifold $M$ and let $L=M/\Delta$ be the (local) quotient space of leaves of $\Delta$.
Then, a vector field $X\in \Delta^\perp$ on $M$ is projectable if and only if
$$
\nabla_TX+C_TX=0\,\,\,\,\mbox{for any}\,\,\,T\in \Delta.
$$
\end{corollary}
\proof For any $T\in \Delta$ we have
$$
[X,T]=(\nabla_XT)_{\Delta}-C_TX-\nabla_TX.
$$
Since $\Delta$ is totally geodesic, then $\nabla_TX\in \Delta^\perp$. Hence $[X,T]\in
\Delta$
if and only if $\nabla_TX+C_TX=0$, and the statement follows from Proposition
\ref{prop:proj}.
\vspace{2ex}\qed

If  $\pi\colon\,M\to L$ is a submersion, we say that
\begin{itemize}
\item[(i)] A one-form $\omega$ on $M$ is {\em projectable\/} if there exist a one-form
$\bar \omega$ on $L$ such that $\bar\omega\circ \pi_*=\omega$.
\item[(ii)] A tensor $D$ on $M$ is {\em projectable\/} if there exists a tensor $\bar D$
on $L$
such that $\bar D\circ \pi_*=\pi_*\circ  D$.
\end{itemize}

 Clearly,  a one-form $\omega$ on $M$ is projectable if and only if  $\omega(X)$ is
constant along the fibers of $\pi$  for any projectable vector field $X$ on $M$. 
Similarly, a tensor $D$ on $M$ is
projectable if and only if $DX$ is projectable for any projectable vector field $X$.

\begin{corollary}\po\label{prop:proj3}
Let $\Delta$ be an  integrable distribution $\Delta$ on
a differentiable manifold $M$,  let $L=M/\Delta$ be the (local) quotient space of leaves
of $\Delta$
and let $\pi\colon\,M\to L$ be the quotient map. Then a one-form $\omega$ on $M$ is
projectable if
and only if $\omega(T)=0$  and $d\omega(T,X)=0$ for any $T\in \Delta$ and  $X\in
\Delta^\perp$.
\end{corollary}
\proof If $\omega=\bar\omega\circ \pi_*$, then $\omega(T)=\bar\omega (\pi_* T)=0$. In
order to
prove that  $d\omega(T,X)=0$
we can assume
that $X$ is projectable. Then $T\omega(X)=0$, hence
$$
d\omega(T,X)=T\omega(X)-X\omega(T)-\omega([X,T])=0
$$
where the vanishing of the last term follows from Proposition \ref{prop:proj}.

 Conversely, if
$X\in \Delta^\perp$ is projectable then
$[X,T]\in \Delta$ by Proposition \ref{prop:proj}, hence the assumptions  give
$$
T\omega(X)=d\omega(T,X)+X\omega(T)+\omega([X,T])=0. \qed
$$

\begin{corollary}\po\label{cor:proj3}
Let $\Delta$ be a totally geodesic distribution on
a Riemannian  manifold $M$ and let $L=M/\Delta$ be the (local) quotient space of leaves
of $\Delta$.
Then a tensor field $D$ on $M$ is projectable if and only if
$$
\nabla_TD=[D,C_T]\,\,\,\,\mbox{for any}\,\,\,\,T\in \Delta.
$$
\end{corollary}
\proof We have
\be\label{eq:condD}
\nabla_TDX+C_TDX=\left( \nabla_TD-[D,C_T]\right)X+D(\nabla_TX+C_TX).
\ee
If $D$ and $X$ are projectable then  $DX$ is also projectable.
Thus, the preceding equality and Corollary \ref{cor:proj} show that $\nabla_TD-[D,C_T]$
vanishes on projectable vector fields, hence it vanishes for this is a tensorial property.

Conversely, if $\nabla_TD=[D,C_T]$, then (\ref{eq:condD}) and Corollary \ref{cor:proj}
imply that $DX$ is projectable whenever $X$ is projectable. \qed

\section[Proof of Proposition \ref{le:1}]{Proof of Proposition \ref{le:1}}

  Let $A^g_\xi$ denote the shape operator of $g$ with respect to  $\xi\in T_g^\perp M$
given by
$$
\<A^g_\xi X, Y\>=\<\alpha_g(X,Y), \xi\>.
$$
Denote by  $D_\xi\colon\,\Delta^\perp\to\Delta^\perp$ the  endomorphism
defined by
$$
D_\xi=A^{-1}A^g_\xi,
$$
where $A$ and $A_\xi^g$ are regarded as endomorphisms of $\Delta^\perp$.

\begin{lemma}\po\label{le:W}
The subspace  of endomorphisms
$W=\spa\{D_\xi\,:\,  \xi\in T_g^\perp M\}$
has dimension two on an open dense subset of  $M^n$.
\end{lemma}

\proof  If $W$ has dimension one at a point, then the first normal space $N_1^g$ of $g$
is also one-dimensional. If this happens on an open subset
$V\subset M^n$ and  $N_1^g$ is not parallel in the normal connection along $V$, then
$V$ is flat by Theorem $1$ in \cite{dt2}. But this contradicts the fact that $f$
has rank two. In case $N_1^g$ is parallel along $V$, then $g(V)$ is contained in a
hyperplane of $\R^{n+2}$, again a contradiction with our assumption that $g$ is
locally substantial.\qed

\begin{lemma}\po\label{le:propD}
The following holds:
\begin{enumerate}
\item[(i)]
$[D_\xi,C_T] = 0\;\;\mbox {for all}\,\,\, T\in\Delta$,
\item[(ii)]
$\nabla_TD_\xi = 0\;\;\mbox {for all}\,\,\, T\in\Delta$ and
$\xi\in T_g^\perp M$ parallel along $\Delta$.
\end{enumerate}
\end{lemma}
\proof
We obtain  from the Codazzi equation that
\be\label{eq:nta}
\nabla_T A=AC_T.
\ee
Moreover,
$$
\nabla_T AD_\xi=\nabla_T A_\xi^g=A_\xi^gC_T=AD_\xi C_T
$$
if $\xi\in T_g^\perp M$ is parallel along $\Delta$. An easy computation yields
\be\label{eq:codtx}
\nabla_T AD_\xi-AD_\xi C_T=(\nabla_TA-AC_T)D_\xi+A(\nabla_TD_\xi-[D_\xi,C_T]),
\ee
hence
\be\label{eq:x}
\nabla_TD_\xi=[D_\xi,C_T].
\ee
On the other hand, we obtain from (\ref{eq:nta}) that $A C_T$ is symmetric, i.e.,
$$
A C_T = C_T^{t} A.
$$
A similar equation holds for $A^g_\xi=A D_\xi$, thus
$$
AD_\xi C_T = A_\xi^gC_T=C_T^{t}A_\xi^g = C_T^{t}AD_\xi =AC_TD_\xi.
$$
This gives $(i)$, and then $(ii)$ in view of (\ref{eq:x}). \vspace{2ex}\qed

We now determine the structure of the splitting tensor $C$. We make use of the
following well-known fact (cf. \cite{dft}).

\begin{proposition}\po\label{le:C=0}
Let $f\colon\,M^n\to \R^{n+1}$ be a hypersurface
of rank two. If the image of the splitting tensor $C$ is either trivial or spanned by the
identity
tensor $I$, then $f$ is surface-like.
\end{proposition}

\begin{remark}\po\label{re:surfacelike} {\em Proposition \ref{le:C=0}
 holds for submanifolds of any codimension. Therefore, since the
splitting tensor $C$ is intrinsic, when $f$ is surface-like its genuine
deformations in $\R^{n+2}$ are surface-like submanifolds over
genuine deformations of $L^2$ in $\R^4$ or $\Sf^4$, respectively.}
\end{remark}

\begin{lemma}\po \label{cor:ct2}
There exists a tensor $J$ on $\Delta^\perp$  such that $J^2=\epsilon I$, $\epsilon\in
\{-1, 1, 0\}$
and  $\spa\,\{I\}\subset C(\Delta)\subset  \spa\,\{I, J\}=W.$
\end{lemma}
\proof  By Proposition \ref{le:C=0} and our assumption that $f$ is not surface-like, on
an open dense subset
 we have that  $C(\Delta)$ is not spanned by $I$. By part $(i)$ of Lemma \ref{le:propD}, 
it is contained in the
subspace $S$ of linear operators on $\Delta^\perp$ that commute with all elements of the 
subspace $W$.
Since $W$ is two-dimensional by Lemma \ref{le:W},  it must  contain   $I$. Otherwise,  it
is easily seen
that $S$ would have to be the subspace spanned by $I$, in contradiction with the fact
that $S$ contains
$C(\Delta)$.
Therefore, $W=\spa\{I, J\}$, where $J$ is a  tensor on $\Delta^\perp$ satisfying
$J^2=\epsilon I$,
$\epsilon\in \{-1, 1, 0\}$.
 In particular, $W\subset S$
and, on the other hand,  the fact that any element of $S$ commutes with  $J$ implies that
the dimension of
$S$ is at most two. Hence $W=S$ and $C(\Delta)\subset S=\spa\,\{I, J\}$.\qed\vspace{1ex}

\begin{lemma}\po \label{le:frame}
There exists a unique (up to signs and permutation) orthonormal frame $\xi_1, \xi_2$ of
$T_g^\perp M$
such that $D_i:=D_{\xi_i}$, $1\leq i\leq 2$, satisfy
$$
\det D_1=1/2=\det D_2.
$$
Moreover,
$D_2^2\neq -D_1^2$ and $\xi_i$, $1\leq i\leq 2$, is parallel along $\Delta$.
\end{lemma}

\proof For any orthonormal frame $\{\xi_1, \xi_2\}$ of $T_g^\perp M$, the Gauss equations
of $f$ and $g$ give
$$
\det D_1 + \det D_2=1.
$$
Since  $W$  has dimension two, we  have $D_1\neq \pm D_2$.

We now show that $D_2^2\neq -D_1^2$. It is easily seen that $D_2^2=-D_1^2$ could only happen if
$W=\spa\{I, J\}$ with $J^2=-I$. Assume  this to be the case.  Let
$\hat
D_i$ denote the complex linear extension of $D_i$ to
$\Delta^\perp\otimes \C$. Then, there exists $\theta\in \Sf^1$ such that
$$
\sqrt{2}\,\hat D_1 = \left[{\begin{array}{cc}
       \theta  &  0  \\
            0     &  \bar\theta
\end{array}} \right]  \,\,\,\,\,\,\,\mbox{and}\,\,\,\,\,\,\,
 \sqrt{2}\,\hat D_2 = \left[{\begin{array}{cc}
       i\theta  &  0  \\
            0     &  -i\bar\theta
\end{array}} \right],
$$
with respect to the frame of $\Delta^\perp\otimes \C$ of eigenvectors of $\hat D_i$,
 $1\leq i\leq 2$. Writing $\theta=e^{i\beta}$ we obtain
$$
\cos\beta D_1-\sin\beta D_2=I\,\,\,\,\,\mbox{and}\,\,\,\,\,\sin\beta D_1+\cos\beta D_2=J.
$$
Hence, the orthonormal frame $\{\xi, \eta\}$ of $T_g^\perp M$ given by
$$
\xi=\cos\beta\, \xi_1-\sin\beta \,\xi_2\,\,\,\,\,\mbox{and}\,\,\,\,\,\eta
=\sin\beta \,\xi_1+\cos\beta\, \xi_2
$$
satisfies
$$
\sqrt{2} A^g_\xi=A\,\,\,\,\,\mbox{and}\,\,\,\,\,\sqrt{2} A^g_\eta=AJ.
$$
Using  the preceding equations, and comparing the Codazzi equation of $f$ with the
Codazzi equation
$$
(\nabla_{X}A^g_\xi)Y-(\nabla_{Y}A^g_\xi)X=\psi(X)A^g_\eta Y-\psi(Y)A^g_\eta X,
$$
where $\psi(Z)=\<\nabla^\perp_Z\xi, \eta\>$, we obtain that $\psi$
vanishes identically. By the Ricci equation, this implies that $[A^g_\xi,
A^g_\eta]=0$, or equivalently, that $[A,AJ]=0$. The latter equation can
only be satisfied if $A=\lambda I$ for some $\lambda\in C^{\infty}(M)$.
 From $\nabla_T A = AC_T$ for any $T\in\Delta$, it follows that
$C_T=\< \grad \lambda, T\> I$. In view of Proposition \ref{le:C=0}, this
is in contradiction with the assumption that $f$ is nowhere surface-like.
Hence $D_2^2\neq -D_1^2$.

Using the preceding condition,  it is easily seen that there exists a unique (up to
signs and permutation) pointwise choice of unit orthogonal  vectors  $\xi_1$ and $\xi_2$ such that $\det
D_1=1/2=\det D_2$,
thus defining a smooth orthonormal normal frame with this property. We now show that
$\xi_1$
and $\xi_2$ are parallel along $\Delta$. Given $x\in M^n$, $T\in \Delta$ and an integral
curve
$\gamma$ of $T$ starting at $x$, let $\hat\xi_i(t)$  denote the parallel transport of
$\xi_i(x)$
along $\gamma$ at $\gamma(t)$. By Lemma \ref{le:propD}-$(ii)$, we have that
$\nabla_{\gamma'(t)}D_{\hat\xi_i(t)}=0$, hence $\det D_{\hat\xi_i(t)}=1/2$.
Since $\xi_1$ and $\xi_2$ are unique (up to signs and permutation) with this property, by
continuity we must have $\hat \xi_i(t)=\xi_i(\gamma(t))$ for any $t$. It
follows that $\nabla_T^\perp\xi_i=0$ for any $T\in \Delta$, $1\leq i\leq 2$. \qed

\begin{lemma}\po\label{ruleparabolic}
There is no open subset $U\subset M^n$ where $W=\spa\,\{I, J\}$
with $J^2=0$. Therefore, the hypersurface $f$ is either elliptic or hyperbolic.
\end{lemma}

\proof It suffices to show  that $f|_U$ is ruled for $U$ as in the statement,
and the proof follows from our assumption.
Let $\{X,Y\}$ be
an orthonormal frame   of  $\Delta^{\perp}$ such that $Y$ spans the image of $J$, that
is, $JY=0$
and $JX=\lambda Y$ for some $\lambda \neq 0$. Since $C(\Delta) \subset W$, we have

\be\label{eq:ctyx}
\<C_TY,X\>=0\;\;\mbox{for all}\,\,T\in\Delta.
\ee
 It follows easily from the fact that $AC_T$ is symmetric that
\be\label{eq:frul}
\<AY,Y\> = 0.
\ee

  We claim that the distribution
$x\mapsto \spa\{Y(x)\}\oplus\Delta(x)$
is totally geodesic.  From (\ref{eq:ctyx}) we have
\be\label{eq:reg1}
\<\nabla_YT,X\> = -\<C_TY,X\>=0,\;\;\mbox{for all}\,\,T\in\Delta.
\ee
Now, let $\{\xi_1, \xi_2\}$ be the orthonormal normal frame given by Lemma \ref{le:frame}.
 From $D_i=D_{\xi_i}\in\spa\{I,J\}$ and $\det D_i=1/2$ we obtain $\sqrt{2}D_iY=Y$, $1\leq i\leq 2$, after replacing $\xi_i$ by $-\xi_i$ if necessary.
Moreover, since $\sqrt{2}D_1$ and $\sqrt{2}D_2$ cannot be both multiples of $I$, we can assume that  $\sqrt{2}D_1X\neq X$.
On the other hand, it follows from Lemma \ref{le:propD}-$(ii)$ that
$\sqrt{2}D_1\nabla_TY=\nabla_T Y$
for all $T\in \Delta$, hence
\be\label{eq:reg2}
\nabla_TY=0\;\;\mbox{for all}\,\,T\in\Delta.
\ee
Now write
$$
A = \left[{\begin{array}{cc}
       \lambda  &  \mu  \\
       \mu  &  0
\end{array}} \right]\;\;\;\;
\mbox{and}\;\;\;\sqrt{2}A_{\xi_1}^g = \left[{\begin{array}{cc}
       \lambda+\theta  &  \mu  \\
       \mu  &  0
\end{array}} \right]
$$
with respect to the frame $\{X,Y\}$ for smooth functions $\lambda, \mu$ and $\theta\neq
0$ on $M^n$.
The Codazzi equations for $A$ and $A^g_{\xi_1}$ yield
$$
Y(\mu)-\lambda\<\nabla_YY,X\>-2\mu\<\nabla_XX,Y\>=0
$$
and
$$
Y(\mu)-(\lambda+\theta)\<\nabla_YY,X\>-2\mu\<\nabla_XX,Y\>=0,
$$
respectively. It follows that
\be \label{eq:reg3}
\<\nabla_YY,X\>=0.
\ee
The claim follows from (\ref{eq:reg1}), (\ref{eq:reg2}), (\ref{eq:reg3}) and the fact that
$\Delta$ is totally geodesic. In view of (\ref{eq:frul}), this implies that $f|_U$ is
ruled.
\vspace{2ex}\qed

 To complete the proof of the direct statement of Proposition \ref{le:1}, 
for the orthonormal frame $\{\xi_1, \xi_2\}$ given by Lemma \ref{le:frame} set
$\phi(Z)=\<\nabla^\perp_Z\xi_1, \xi_2\>$.
Then, the last assertion in Lemma \ref{le:frame} shows that condition $(i)$ holds, and
hence $(iv)$ in view of Lemma \ref{le:propD}.
Condition $(v)$ follows from the Codazzi equation for $g$, whereas
$(vi)$ and  $(vii)$   from the  Ricci equation.

   The proof of the converse statement   is a straightforward application of the
fundamental
theorem of submanifolds. Choose an orthonormal frame $\{\xi_1, \xi_2\}$ of the trivial
bundle
$E=M\times \R^2$ and define a connection $\hat\nabla$ on $E$ by
$\<\hat\nabla_X\xi_1, \xi_2\>=\phi(X)$ for $X\in TM$. It follows from $(i)$  that
$\xi_1$ and $\xi_2$ are parallel along $\Delta$. Let $\alpha\colon\,TM\times TM\to E$ be defined
by setting $\ker\alpha=\Delta$ and
$$
\alpha(X,Y)=\sum_{i=1}^2\<AD_iX, Y\>\xi_i\,\,\,\mbox{for all}\,\,\,\, X, Y\in
\Delta^\perp.
$$
Condition $(ii)$ implies that $\alpha$ is symmetric, and from $(iii)$ and the Gauss
equation
for $f$  it satisfies the Gauss equation for an isometric immersion  into $\R^{n+2}$.
The Codazzi equation follows from $(iv)$,  $(v)$ and (\ref{eq:codtx}), whereas the Ricci
equation
is a consequence of $(vi)$ and $(vii)$. By the fundamental theorem of submanifolds, there
exists
an isometric immersion $g\colon\,M^n\to \R^{n+2}$ having $\alpha$ as second fundamental
form and
$\hat\nabla$ as normal connection. Since $D_1\neq \pm D_2$, it follows that the first
normal
spaces of $g$ have dimension $2$ everywhere, hence $g$ is locally substantial.

 Finally, the last assertion is a consequence of the uniqueness
of  the
frame  $\{\xi_1,\xi_2\}$ such that $\det D_1=1/2=\det D_2$,
together with the uniqueness part of the fundamental theorem of submanifolds.

\section[Proof of Proposition \ref{le:2}]{Proof of Proposition \ref{le:2}}

We start with the following lemma.

\begin{lemma}\po\label{le:phi}
The tensors $D_1$, $D_2$ and the one-form $\phi$ are
projectable with respect to the quotient map $\pi\colon\,M^n\to L^2=M/\Delta$.
\end{lemma}
\proof That $D_1$ and $D_2$ are  projectable  follows from  $(iv)$ and Corollary
\ref{cor:proj3}.
On the other hand, $\phi$ is projectable by $(vi)$ and Corollary
\ref{prop:proj3}.\vspace{2ex}\qed

Hence, there exist tensors $\bar D_i$, $1\leq i\leq 2$, and a one-form
$\bar \phi$ on $L^2$ such that $\bar D_i\circ \pi_*=\pi_*\circ D_i$ and
$\bar \phi\circ \pi_*=\phi$. In particular, $[\bar D_1, \bar D_2]\circ \pi_*=\pi_*\circ
[D_1,  D_2]=0$,
hence there exists a unique tensor $\bar J$ on $L^2$ such that  $\bar J^2=\e I$, $\e\in
\{1, -1\}$,
and $\bar D_i\in \{I, \bar J\}$. Write $\bar J=a\bar D_1 + b\bar D_2$, $a, b\in \R$, and
define
$\hat J=aD_1 +bD_2$. From $\bar D_i\circ \pi_*=\pi_*\circ D_i$ we obtain that
$\hat J\circ \pi_*=\pi_*\circ \bar J$, hence $\hat J^2=\epsilon I$,
$\epsilon \in \{1, -1\}$. Since $J$ is (up to sign) the unique tensor in $\spa\{D_1,
D_2\}$
with this property, it follows that $\hat J=J$, after a change of sign if necessary.
Summarizing, we have proved the existence of a unique tensor $\bar J$ on $L^2$
such that  $\hat J\circ \pi_*=\pi_*\circ \bar J$ and $\bar D_i\in \spa\{I, \bar J\}$.

Conditions $(a)$ and $(d)$ are clear, for these properties are inherited from
$D_1$ and $D_2$. In order to verify the remaining conditions we first
make a few computations.

Let $X$  and $Y$ be projectable vector fields on $M^n$. Then, we have
\be\label{eq:p1}
f_*AX=-N_*X=-h_*\pi_*X,
\ee
$$
f_*AD_iX=-h_*\pi_*D_iX=-h_*\bar D_i\pi_*X
$$
and
\begin{eqnarray}\label{eq:p3}
f_*\nabla_{X}A D_iY\!\! &=& \!\!\nabla_{X}  f_* A D_iY - \<A X,AD_i Y\>N\nonumber\\
\!\!&=&\!\! -{\nabla}_{\pi_*X}h_*\bar D_i\pi_*Y - \<h_*\pi_*X,h_*\bar D_i\pi_*Y\>h\circ
\pi\nonumber\\
\!\!&=&\!\! -h_*\nape_{\pi_*X}\bar D_i\pi_*Y- \alpha_h(\pi_*X,\bar D_i\pi_* Y).
\end{eqnarray}
In view of (\ref{eq:p1})--(\ref{eq:p3}), condition   $(v)$  implies $(b)$ and
\be\label{eq:adi}
\alpha_h(\bar D_iX,Y)=\alpha_h(X, \bar D_i Y),\,\,\,1\leq i\leq 2,
\ee
which is equivalent to $h$ being elliptic or hyperbolic with respect to $\bar J$.
On the other hand, condition $(vii)$ gives $(c)$.

To prove (\ref{eq:hessgamma}), let $\hat \pi\colon\,T_h^\perp L\to L$ be the
canonical projection.
 By the Gauss parametrization, there exists a local diffeomorphism
$\Phi\colon\,U\subset T_h^\perp L\to M$ of an open neighborhood of the zero
section such that  $\pi\circ\Phi=\hat \pi$ and
$$
\psi(x,w):=f\circ \Phi(x,w)=\gamma h +h_*\nabla \gamma+w.
$$
For any horizontal vector $X\in T_{(x,w)}(T_h^\perp L)$ we have
$$
\psi_*X=h_*P\hat\pi_*X+\alpha_h(\hat\pi_*X,\nabla \gamma),
$$
 where $P$ is the endomorphism of $TL$ given by
\be\label{eq:P}
P={\rm Hess}_{\,\gamma}+\gamma I-B_w.
\ee
Here $B_w$ stands for the shape operator of $h$ in direction $w$. Thus,
$$
f_*\Phi_*X=h_*P\hat \pi_*X+\alpha_h(\hat\pi_*X,\nabla \gamma)
=h_*P\pi_*\Phi_*X+\alpha_h(\pi_*\Phi_*X,\nabla \gamma),
$$
and hence
\be\label{eq:adsym}
-\<AD_i\Phi_*X,\Phi_*Y\>=\<h_*\bar D_i\pi_*\Phi_*X,h_*P\pi_*\Phi_*Y\>=\<\bar D_i\hat
\pi_*X,P\hat \pi_*Y\>'
\ee
for all horizontal vectors $X, Y\in T(T_h^\perp L)$.
Therefore, condition $(ii)$  implies that $P\bar D_i=\bar D_i^tP$. Since
$B_w\bar D_i=\bar D_i^tB_w$ by $(\ref{eq:adi})$, this gives (\ref{eq:hessgamma}).

 We now  prove the converse.  Set $\omega=\bar\omega\circ \pi_*$ and let $D_i$ be the
horizontal lift
of $\bar D_i$ to $M^n$, $1\leq i\leq 2$, that is, $\Delta\subset \ker D_i$ and, for any
$x\in M^n$ and
$X\in \Delta^\perp(x)$,   $D_iX$ is the unique vector in $\Delta^\perp(x)$ such that
$\pi_*D_iX=\bar D_i\pi_*X$.
Define in a similar way a tensor $J$ on  $\Delta^\perp$ such that $\pi_*\circ J=\bar
J\circ \pi_*$,
so that  $\spa\{D_1, D_2\}=\spa\{I, J\}$.

 Conditions $(i)$,   $(iii)$ and $(viii)$  are clear. Since $B_w\bar J=\bar J^t B_w$ for
any $w\in T_h^\perp L$,
for $h$ is elliptic or hyperbolic with respect to $\bar J$, it follows from
(\ref{eq:hessgamma})
that $P\bar J=\bar J^t P$, where $P$ is given by (\ref{eq:P}). This implies that
$P\bar D_i=\bar D_i^t P$ for $1\leq i\leq 2$, hence $AD_i=D_i^tA$ by (\ref{eq:adsym}).
This proves $(ii)$.

  It follows from Corollary \ref{cor:proj3} that $\nabla_TD_i=[C_T, D_i]$, hence
(\ref{eq:codtx}) gives
$\nabla_TAD_i=AD_iC_T$. This implies that $AD_iC_T$ is symmetric, which is equivalent to
$A[D_i, C_T]=0$,
bearing in mind $(ii)$  and the fact that $AC_T$ is symmetric.
This proves $(iv)$. Moreover, it implies that $C(\Delta)\subset \spa\{I, J\}$, hence $f$
is elliptic
or hyperbolic, according as $J^2=-I$ or $I$, that is, according as $h$ is elliptic or
hyperbolic, respectively.

   Since (\ref{eq:adi}) is satisfied, for $h$ is elliptic or hyperbolic with respect to
$\bar J$,
using (\ref{eq:p1})--(\ref{eq:p3}) and  condition $(b)$ we obtain  $(v)$. Finally,
condition $(vi)$
follows from Corollary \ref{prop:proj3} and $(vii)$ is a consequence of $(c)$, by using
(\ref{eq:p1}).\qed

\section[Proof of Proposition \ref{le:3}]{Proof of Proposition \ref{le:3}}

We make use of the following special case of Theorem $5$ in \cite{dt1}.

\begin{proposition}\po\label{prop:comp}
Let $f\colon\, M^n\to \R^{n+1}$ and $g\colon\, M^n\to\R^{n+2}$
be isometric immersions. If $g$ is the composition $g=H\circ f$ of $f$ with an isometric
immersion
$H\colon\, W\to \R^{n+2}$ of an open subset $W\supset f(M)$ of $\R^{n+1}$, then there
exists an orthonormal
frame $\{\xi, \eta\}$ of $T_g^\perp M$ such that $A_\xi^g=A^f$ and $\rank A_\eta^g\leq 1$.
The converse also holds if $\rank A_\eta^g=1$ everywhere and $\eta$ is parallel along
$\ker A_\eta$.
\end{proposition}

  For $f$ and $g$ as in Proposition \ref{le:1}, assume that there exist an open subset
$U\subset M^n$
and an isometric immersion $H\colon\, W\to \R^{n+2}$ of an open subset $W\supset f(U)$ of
$\R^{n+1}$ such
that $g|_U=H\circ f|_U$. Then, by  Proposition \ref{prop:comp}  there exists $\theta\in
(0, 2\pi)$
such that
$$
D_\theta:=\cos\theta D_1+\sin\theta D_2=I\,\,\,\mbox{and}\,\,
\rank (D_{\theta+\pi/2}:=-\sin\theta D_1+\cos\theta D_2)<2.
$$
Since this can never happen if $D_1, D_2\subset \{I, J\}$ with $J^2=-I$, the last assertion is proved.

  From now on assume that $D_1, D_2\subset \{I, J\}$ with $J^2=I$.
Let $\{X,Y\}$ be a frame of $\Delta^\perp$ of eigenvectors of $J$, say,  
$$
\sqrt{2}\, D_1 = \left[{\begin{array}{cc}
       \theta_1  &  0  \\
            0     &  1/\theta_1
\end{array}} \right]  \,\,\,\,\,\,\,\mbox{and}\,\,\,\,\,\,\,\sqrt{2}\, D_2 =
\left[{\begin{array}{cc}
       \theta_2  &  0  \\
            0     &  1/\theta_2
\end{array}} \right],\,\,\,\,\theta^2\neq \pm\theta^1.
$$
Then, it is easily checked that
$a_1D_1+a_2D_2=I$ if and only if
$$
a_i=\frac{\sqrt{2}\theta_i(1-\theta^2_j)}{\theta_i^2-\theta_j^2}, \,\,1\leq i\neq j\leq 2.
$$
Moreover, for these values of $a_1$ and $a_2$, the rank of $-a_2D_1+aD_2$ is less than
two if and
only if either $\theta_1^2+\theta_2^2=2$ or $\theta_1^{-2}+\theta_2^{-2}=2$, that is, if
and only
if $\rank D_1^2+D_2^2-I<2$. Note that $\sqrt{2} a_i=\theta_i$ if
$\theta_1^2+\theta_2^2=2$, whereas $\sqrt{2} a_i=1/\theta_i$ if
$\theta_1^{-2}+\theta_2^{-2}=2$. In either case we have $a_1^2+a_2^2=1$.

Summarizing, there exists an orthonormal frame $\{\xi, \eta\}$ of $T_g^\perp M$ such that
$A_\xi^g=A^f$ and $\rank A_\eta^g\leq 1$ if and only if $\rank D_1^2+D_2^2-I<2$.
This already shows the ``if" part of the statement.

Let us prove the converse. Assume, say,  that $\theta_1^2+\theta_2^2=2$ on an open subset
$U\subset M^n$.
Then $\sqrt{2}(\theta_1 D_1+\theta_2D_2)=I$ and $X$ belongs to the kernel of $-\theta_2
D_1+\theta_1D_2$.
Therefore, the orthonormal frame $\{\xi, \eta\}$ of $T_g^\perp M$ given by
$$
\sqrt{2}\,\xi=\theta_1 \xi_1+\theta_2\xi_2\,\,\,\mbox{and}\,\,\,\,\sqrt{2}\,\eta
=-\theta_2 \xi_1+\theta_1\xi_2
$$
satisfies $A_\xi^g=A$ and $\rank A_\eta^g=1$.
By Proposition \ref{prop:comp}, in order to conclude that $g|_U$ is a composition
$g|_U=H\circ f|_U$,
where $H\colon\, W\to \R^{n+2}$ is an isometric immersion of an open subset $W\supset
f(U)$,
we must still show that $\eta$ is parallel along $\ker A_\eta^g$, that is, that
$\nabla_X^\perp \eta=0$. The latter condition is equivalent to
\be\label{eq:parallel}
X(\theta_1)=\theta_2\phi(X)
\ee
bearing in mind that $\theta_1^2+\theta_2^2=2$.

In order to prove (\ref{eq:parallel}) it is easier to work on the quotient space $L^2$.
First note that Proposition \ref{le:1}-$(iv)$ implies that $\theta_1$ and $\theta_2$ are
constant along the leaves of $\Delta$, hence give rise to functions on $L^2$ that we also
denote by $\theta_1$ and $\theta_2$. Moreover, these are also the eigenvalues of $\bar
D_1$
and $\bar D_2$. Let $(u,v)$ be coordinates in $L^2$ whose coordinate vector fields are
eigenvectors of $\bar D_i$, i.e.,
\be\label{eq:cs}
\sqrt{2}\,\bar D_1 = \left[{\begin{array}{cc}
       \theta_1  &  0  \\
            0     &  1/\theta_1
\end{array}} \right]  \,\,\,\,\,\,\,\mbox{and}\,\,\,\,\,\,\,\sqrt{2}\,\bar D_2
= \left[{\begin{array}{cc}
       \theta_2  &  0  \\
            0     &  1/\theta_2
\end{array}} \right],
\ee
with respect to  the frame $\{\d_u, \d_v\}$  of coordinate vector fields.
Then, equation (\ref{eq:parallel}) is equivalent to
\be\label{eq:parallel2}
(\theta_1)_u=\phi^u\theta_2.
\ee
The equation in Proposition \ref{le:2}-$(b)$ can be written as
\be\label{eq:cod3a}
\nabla_{\d_u}\bar D_i\d_v-\nabla_{\d_v}\bar D_i\d_u=(-1)^{j}(\bar \phi^u \bar D_j\d_v
-\bar\phi^v \bar D_j\d_u),\,\,\,1\leq i\neq j\leq 2,
\ee
where $\bar\phi^u=\bar\phi(\d_u)$ and  $\bar\phi^v=\bar\phi(\d_v)$. This is equivalent to

$$
\frac{\theta^i_u}{(\theta^i)^2}+\left(\theta^i-\frac{1}{\theta^i}\right)\Gamma^v
=(-1)^{i}\frac{\phi^u}{\theta^{j}}
$$
and
$$
{\theta^i_v}+\left(\theta^i-\frac{1}{\theta^i}\right)\Gamma^u
=(-1)^{j}{\theta^{j}}\phi^v,\,\,\,\,1\leq i\neq j\leq 2,
$$
where   we write $\Gamma^u=\Gamma_{uv}^u$ and $\Gamma^v=\Gamma_{uv}^v$ for simplicity.
In terms of $\tau^i:=(\theta^i)^2$, the preceding equations become
\be\label{eq:hyp1a}
\left(\frac{1}{\tau^i}\right)_u+2\left(\frac{1}{\tau^i}-1\right)\Gamma^v=
2(-1)^{j}\frac{\phi^u}{\theta_1\theta_2}
\ee
and
\be\label{eq:hyp2a}
\tau^i_v+2(\tau^i-1)\Gamma^u=2(-1)^{j}\phi^v\theta_1\theta_2,\,\,\,\,1\leq i\neq j\leq 2.
\ee
Equation (\ref{eq:parallel2}) takes the form
$$
(\tau_1)_u=2\phi^u\theta_1\theta_2,
$$
and we must show that it is satisfied if $\tau_1+\tau_2=2$. We obtain from
(\ref{eq:hyp1a}) for $i=1$  that
\be\label{eq:phiu1}
2\phi^u\theta_1\theta_2=-\frac{\tau_2}{\tau_1}(\tau_1)_u+2(2-\tau_1)(1-\tau_1)\Gamma^v.
\ee
On the other hand, equation (\ref{eq:hyp1a}) for $i=2$ gives
\be\label{eq:phiu2}
2\phi^u\theta_1\theta_2=-\frac{\tau_1}{\tau_2}(\tau_1)_u+2\tau_1(1-\tau_1)\Gamma^v.
\ee
It follows from (\ref{eq:phiu1}) and (\ref{eq:phiu2}) that
$$
\frac{(\tau_1)_u}{\tau_2}=\tau_1(1-\tau_2)\Gamma^v.
$$
Replacing into (\ref{eq:phiu2}) yields
$$
2\phi^u\theta_1\theta_2=-\frac{\tau_1}{\tau_2}(\tau_1)_u+2\frac{(\tau_1)_u}{\tau_2}
=\frac{2-\tau_2}{\tau_1}(\tau_1)_u=(\tau_1)_u.\,\,\,\qed
$$

\section[Proof of Proposition \ref{le:4}]{Proof of Proposition \ref{le:4}}

In order to prove the direct statement, we consider separately the hyperbolic and
elliptic cases.

\subsection[The hyperbolic case]{The hyperbolic case}

As in the proof of Proposition \ref{le:3}, let  $(u,v)$ be conjugate coordinates on $L^2$
such that $\bar D_1$ and $\bar D_2$ are given by (\ref{eq:cs})
with respect to  the frame $\{\d_u, \d_v\}$  of coordinate vector fields.   As pointed out
in the proof of Proposition \ref{le:3},   condition $(b)$ can be written as
(\ref{eq:hyp1a}) and (\ref{eq:hyp2a}) in terms of $\tau^i:=(\theta^i)^2$, $1\leq i\leq 2$.
On the other hand, $(c)$ takes the form
\be\label{eq:hyp3a}
2(\phi^v_u-\phi^u_v)= \frac{\tau^1-\tau^2}{\theta^1\theta^2}F.
\ee
It follows from (\ref{eq:hyp1a}) and  (\ref{eq:hyp2a}) that
\be\label{eq:hyp1b}
\left(\frac{1}{\tau^1}+ \frac{1}{\tau^2}\right)_u+2\left(\frac{1}{\tau^1}+
\frac{1}{\tau^2}-2\right)\Gamma^v=0
\ee
and
\be\label{eq:hyp2b}
(\tau^1+\tau^2)_v+2(\tau^1+\tau^2-2)\Gamma^u=0.
\ee
In terms of
$$
\alpha=\tau_1+\tau_2\,\,\,\mbox{and}\,\,\,\beta=1/\tau_1+1/\tau_2
$$
the preceding equations can be written as
\be\label{eq:ab}
\beta_u + 2(\beta-2)\Gamma^v=0 \,\,\,\mbox{and}\,\,\,\, \a_v + 2(\a-2)\Gamma^u=0.
\ee
Notice that $\alpha,\beta>0$. Moreover, since $\tau_1$ and $\tau_2$ are distinct
real roots of $\tau^2-\alpha\tau+(\alpha/\beta)=0$, it follows that $\alpha\beta>4$
and that  $\tau_1$ and $\tau_2$ can be
recovered from $\alpha$ and $\beta$ by
\be\label{eq:recover}
2\tau_i=\a-(-1)^i\sqrt{(\a/\beta)(\a\beta-4)}, \,\,\,\,\,1\leq i\leq 2.
\ee
Since $\bar D_1$ and $\bar D_2$ satisfy the condition in Proposition \ref{le:3},
we have that $\alpha\neq 2$ and $\beta\neq 2$. Then, we can define
\be\label{eq:gh}
\varphi=1/(\a-2)\,\,\,\mbox{and}\,\,\,\, \psi=1/(\beta-2).
\ee
 From $\alpha>0$, $ \beta>0$ and $ \alpha\beta-4>0$, it follows that $(\varphi, \psi)$
satisfies
(\ref{eq:cond1}). Moreover, we obtain from (\ref{eq:ab}) that
$$
\varphi_v/\varphi=2\Gamma^u\,\,\,\mbox{and}\,\,\,\, \psi_u/\psi=2\Gamma^v.
$$
Set
$$
\rho:=\sqrt{|2(\varphi+\psi)+1|}={\sqrt{\a\beta-4}}/{\sqrt{|(\a-2)(\beta-2)|}}.
$$
Writing $\phi^u$, $\phi^v$ and the $\tau^i$ in terms of $\alpha$ and $\beta$   by means of
(\ref{eq:hyp1a}), (\ref{eq:hyp2a}) and (\ref{eq:recover}), and replacing into
(\ref{eq:hyp3a}),
a rather long but straightforward computation shows that it reduces to $Q(\rho)=0$.
Thus, the set ${\cal C}_h$ is nonempty. Moreover, (\ref{eq:hessgamma})  reduces to
$$
{\rm Hess}_{\,\gamma}(\d_u, \d_v)+F\gamma=0,
$$
that is, to $Q(\gamma)=0$.  Finally, distinct triples $(\bar D_1, \bar D_2, \bar \phi)$
(up to signs and permutation of $\bar D_1$ and $\bar D_2$) give rise to distinct
$4$-tuples
$(\tau^1, \tau^2, \phi^u, \phi^v)$, and hence to distinct pairs $(\varphi, \psi)$.

\subsection[The elliptic case]{The elliptic case}

Assume that    $(u,v)$ are coordinates on $L^2$ such that the frame $\{\d_z, \d_{\bar
z}\}$
defined  by
$$
\d_z=\frac{1}{2}(\d_u-i\d_v)\,\,\,\mbox{and}\,\,\,\,\d_{\bar z}=\frac{1}{2}(\d_u+i\d_v),
$$
in  terms of the   frame $\{\d_u, \d_v\}$ of coordinate vector fields, are eigenvectors of
the complex linear extension of $\bar J$ to  $TL\otimes \C$. Write the complex linear
extensions of $\bar D_i$ as
\be\label{eq:cs2}
\sqrt{2}\,\bar D_1 = \left[{\begin{array}{cc}
       \theta^1  &  0  \\
            0     &  \bar\theta^1
\end{array}} \right],  \,\,\,\,\,\,\, \sqrt{2}\,\bar D_2 = \left[{\begin{array}{cc}
       \theta^2  &  0  \\
            0     &  \bar\theta^2
\end{array}} \right],\,\,\,\,|\theta^1|=1=|\theta^2|, \,\,\theta^2\neq \pm\theta^1,
\ee
with respect to  the frame $\{\d_z, \d_{\bar z}\}$.
Define a complex valued Christoffel symbol $\Gamma$ by
$$
\nabla_{\d_ z}\d_{\bar z}=\Gamma \d_z+ \bar \Gamma \d_{\bar z}
$$
and set $\phi^z=\phi(\d_z)$. As in the hyperbolic case, define $\tau^i=(\theta^i)^2$,
$1\leq i\leq 2$.
Then, the complex versions of  (\ref{eq:hyp1a}) and (\ref{eq:hyp2a}) are
\be\label{eq:444}
\tau^i_{\bar z}+2(\tau^i-1)\Gamma=2(-1)^j\bar\phi^z\theta^1\theta^2,
\ee
whereas (\ref{eq:hyp3a}) becomes
\be\label{eq:hyp4a}
4\mbox{Im}\,(\phi^z)_{\bar z}=i(\tau^1-\tau^2)\bar \theta^1\bar \theta^2 F.
\ee
Notice that, since $\bar \theta^i=1/\theta^i$ these are the same equations
as in the hyperbolic case with $(u,v)$,   $(\phi^u, \phi^v)$ and $(\Gamma^u, \Gamma^v)$
replaced by $(z, \bar z)$,  $(\phi^z, \bar\phi^z)$ and  $(\Gamma, \bar \Gamma)$,
respectively.
Set $\alpha=\tau^1+\tau^2$ as before. We obtain from (\ref{eq:444}) that
$$
\alpha_{\bar z}+2(\alpha-2)\Gamma=0.
$$
Note that $|\alpha|<2$, since
$|\tau^1|=1=|\tau^2|$ and $\tau^1\neq \tau^2$.
 Thus $\varphi=1/(\a-2)$ is
well-defined and satisfies
$$
\frac{\varphi_{\bar z}}{\varphi}=2\Gamma.
$$
   From $|\alpha|<2$  and
$$
4\mbox{Re}(\varphi) +1=\frac{|\alpha|^2-4}{|\alpha-2|^2}
$$
it follows that $4\mbox{Re}(\varphi) +1< 0$. On the other hand, condition $(d)$ implies
that
$\tau^2\neq -\tau^1$,  hence $\alpha\neq 0$. This gives $\varphi\neq -1/2$, thus
conditions (\ref{eq:cond2}) are satisfied. Observe that, since $\alpha\neq 0$, we can
recover
$\tau^1$ and $\tau^2$ from  $\alpha$ by
\be\label{eq:rec2}
\tau^i=\frac{\alpha}{2}\left(1\pm i\frac{\sqrt{4-|\alpha|^2}}{|\alpha|}\right).
\ee
Now set
$$
\rho=\sqrt{-(4\mbox{Re}(\varphi) +1)}=\frac{\sqrt{4-|\alpha|^2}}{|\alpha-2|}.
$$
As before, writing $\phi^z$ and the $\tau^i$ in terms of $\alpha$ by means
of (\ref{eq:444}) and (\ref{eq:rec2}), and replacing into (\ref{eq:hyp4a}),
we arrive at the equation  $Q(\rho)=0$. This shows that the set ${\cal C}_h$
is nonempty. Finally, (\ref{eq:hessgamma})  reduces in this case to
$$
{\rm Hess}_{\,\gamma}(\d_z, \d_{\bar z})+F\gamma=0,
$$
that is, to $Q(\gamma)=0$. Again, distinct triples $(\bar D_1, \bar D_2, \bar \phi)$
(up to signs and permutation of $\bar D_1$ and $\bar D_2$) yield distinct
triples $(\tau^1, \tau^2, \phi^z)$,  and hence distinct  $\varphi$'s.

\subsection[Proof of the converse statement]{Proof of the converse statement}

 We argue first in the hyperbolic case.
Let $(\varphi, \psi)$ be a pair of smooth functions on $L^2$ satisfying (\ref{eq:vp}) and
 (\ref{eq:cond1}).
Assume also that   $\rho=\sqrt{|2(\varphi+\psi)+1|}$ satisfies $Q(\rho)=0$.  Set
$$
\alpha=2+1/\varphi\,\,\,\mbox{and}\,\,\, \beta=2+1/\psi.
$$
Since $(\varphi, \psi)$ satisfies  (\ref{eq:cond1}), it follows that $\alpha>0$,
$\beta>0$ and $\alpha\beta-4>0$.
Then, we can define  $\tau^i$  by
(\ref{eq:recover}).   Let $\phi^u$ and $\phi^v$ be given by  (\ref{eq:hyp1a}) and
(\ref{eq:hyp2a}), respectively.
 It follows from $Q(\rho)=0$ that (\ref{eq:hyp3a}) is satisfied.  Write
$\tau^i=(\theta^i)^2$, 
let  $\bar D_1$ and $\bar D_2$ be defined by  (\ref{eq:cs}) with respect to  the frame
$\{\d_u, \d_v\}$
of coordinate vector fields, and set $\bar \phi=\phi^u\, d u+ \phi^v\, d v$. Then
condition $(a)$ is clear,
whereas $(b)$ follows from (\ref{eq:hyp1a}) and (\ref{eq:hyp2a}). Condition $(c)$ is a
consequence
of  (\ref{eq:hyp3a}),  (\ref{eq:hessgamma}) follows from  $Q(\gamma)=0$ and $(d)$ is
automatic  in this case.
It is also clear that distinct pairs $(\varphi, \psi)$ give rise to distinct 4-tuples
$(\tau^1, \tau^2,\phi^u, \phi^v)$, and hence to distinct triples $(\bar D_1, \bar D_2,
\bar \phi)$.

   Assume now that  $h$ elliptic with  complex conjugate coordinates $(z, \bar z)$.
Suppose that $\varphi$ satisfies (\ref{eq:cond2}) and   
$\rho:=\sqrt{-(4\mbox{Re}(\varphi) +1)}$
satisfies $Q(\rho)=~0$.  Set
$$
\alpha=2+1/\varphi.
$$
It follows from (\ref{eq:cond2}) that $\alpha\neq 0$ and $|\alpha|<2$.
Then (\ref{eq:rec2}) gives   $\tau^i$,  $1\leq i\leq 2$,  with $\tau^2\neq \pm\tau^1$,
$|\tau^1|=1=|\tau^2|$ and  $\tau^1+\tau^2=\alpha$.
Write $\tau^i=(\theta^i)^2$ and define a triple  $(\bar D_1, \bar D_2, \bar \phi)$
by requiring that  $\bar D_1$ and $\bar D_2$ be given by (\ref{eq:cs2}) with respect to
the
frame $\{\d_z, \d_{\bar z}\}$,  and $\bar \phi=\phi^z\, d z+ \bar\phi^z\, d \bar z$ where
$\phi^z$
is given by (\ref{eq:444}). Then, condition $(a)$
is immediate and $(b)$ follows from (\ref{eq:444}). From  $Q(\rho)=0$ we obtain
(\ref{eq:hyp4a}), and hence condition
$(c)$. Finally, (\ref{eq:hessgamma}) follows from $Q(\gamma)=0$, and $(d)$ holds because
$\tau^2\neq -\tau^1$.
Again, distinct $\varphi$'s yield distinct triples $(\bar D_1, \bar D_2, \bar \phi)$.\qed

{\renewcommand{\baselinestretch}{1}
\hspace*{-30ex}\begin{tabbing}
\indent \= IMPA  \hspace{30ex} Universidade Federal de S\~ao Carlos \\
\>  Estrada Dona Castorina, 110 \hspace{7ex}
Via Washington Luiz km 235 \\
\> 22460-320 --- Rio de Janeiro
\hspace{8ex} 13565-905 --- S\~ao Carlos  \\
\> Brazil\hspace{31ex} Brazil\\
\> marcos@impa.br,\,\, luis@impa.br  \hspace{5ex}
tojeiro@dm.ufscar.br
\end{tabbing}}

\end{document}